\documentclass[11pt]{article}

\usepackage{color}
 \usepackage[T1]{fontenc}
\usepackage[utf8]{inputenc}

\newtheorem{myproposition}{Proposition}[section]
\newtheorem{mytheorem}[myproposition]{Theorem}
\newtheorem{mylemma}[myproposition]{Lemma}

\newtheorem{myconjecture}[myproposition]{Conjecture}

\newtheorem{myfact}[myproposition]{Fact}

\newtheorem{myproblem}[myproposition]{Problem}
\usepackage{amsfonts}
\usepackage{amsmath}
\usepackage{authblk}

\newcommand{\qed}{\hfill \rule{.1in}{.1in}}
\date{}
\makeatletter
\def\imod#1{\allowbreak\mkern10mu({\operator@font mod}\,\,#1)}
\makeatother

\begin {document}

\title{Total vertex product irregularity strength of graphs}
\author[1]{Marcin Anholcer}
\author[2]{Azam Sadat Emadi}
\author[2]{Doost Ali Mojdeh}
\affil[1]{\scriptsize{}Pozna\'n University of Economics and Business, Institute of Informatics and Quantitative Economics}
\affil[ ]{Al. Niepodleg{\l}o\'sci 10, 61-875 Pozna\'n, Poland, \textit{m.anholcer@ue.poznan.pl}}
\affil[2]{Department of Mathematics, University of Mazandaran, Babolsar, Iran\\
\textit{math\_emadi2000@yahoo.com}, \textit{damojdeh@umz.ac.ir}}
\maketitle
\begin{abstract}
Consider a simple graph $G$. We call a labeling $w:E(G)\cup V(G)\rightarrow \{1, 2, \dots, s\}$ (\textit{total vertex}) \textit{product-irregular}, if all product degrees $pd_G(v)$ induced by this labeling are distinct, where $pd_G(v)=w(v)\times\prod_{e\ni v}w(e)$. The strength of $w$ is $s$, the maximum number used to label the members of $E(G)\cup V(G)$. The minimum value of $s$ that allows some irregular labeling is called \textit{the total vertex product irregularity strength} and denoted $tvps(G)$. We provide some general bounds, as well as exact values for chosen families of graphs.
\end{abstract}

\section{Introduction}

Let $G$ be a simple undirected graph  with no loops with vertex set $V(G)$, edge set $E(G)$, minimum degree $\delta(G)$ and maximum degree $\Delta(G)$. Let us assign a label (positive integer) $w(x)$ to every $x\in E(G)\cup V(G)$. For each $v \in V(G)$, its \textit{product degree} will be defined by the formula

$$
pd_G(v)=w(v)\times \prod_{e\ni v}w(e),
$$

\noindent{}where $d_G(v)$ is the degree of $v$ in $G$ and $e\ni v$ means that the vertex $v$ is incident to the edge $e$. In particular, in the case of an isolated vertex $v$ we have $pd_G(v)=w(v)$.

We call $w$ (\textit{total vertex}) \textit{product-irregular} when for each couple of vertices $u,v\in V(G), u\neq v$, we have $pd_G(u)\neq pd_G(v)$. The (\textit{total vertex product}) \textit{strength} of the labeling $w$ is defined as

$$
tvps_w(G)=\max\{w(x)|x\in E(G)\cup V(G)\}.
$$
Consequently, the \textit{total vertex product irregularity strength} of $G$ is described with the formula
$$
tvps(G)=\min\{tvps_w(G)|w\, \mbox{\rm{is product-irregular}}\}.
$$

This concept is a variant of the product irregularity strength, introduced by Anholcer \cite{ref_Anh1} and studied in \cite{ref_Anh2} and \cite{ref_DarHuj}. This time the graphs under consideration have no isolated edges and at most one isolated vertex and the product degree is defined as the product of the labels of the incident edges only:
$$
pd_G(v)=
\begin{cases}
\prod_{e\ni v}w(e) \text{     if     }d_G(v)>0,\\
0\text{     if     }d_G(v)=0.
\end{cases}
$$
Similarly as in the definition above, one defines the (product) strength of the labeling $ps_w(G)$ as the maximum label used and the respective graph invariant (\textit{product irregularity strength}, $ps(G)$) is the minimum strength among all the product-irregular labelings.

The motivation to study this kind of problems were the well-known \textit{irregularity strength} and \textit{total vertex irregularity strength}, where the vertex weighted degrees are defined as the sums of labels (instead of products). The irregularity strength was defined in \cite{ref_ChaJacLehOelRuiSab1} and then analyzed by many authors (see e.g. \cite{ref_AigTri,ref_AmaTog,ref_AnhPal,
ref_FauLeh1,ref_FerGouKarPfe,ref_FriGouKar1,ref_Leh,ref_Nie1}). The best known general upper bound is $s(G)\leq 6n/\delta$ for graphs with $\delta\geq 6$ (see \cite{ref_KalKarPfe1}). Majerski and Przyby{\l}o \cite{ref_Prz3} improved it for dense graphs of sufficiently large order (they proved that in such case $s(G)\leq (4+o(1))n/\delta+4$). The total vertex irregularity strength was in turn introduced by Baca et al. in \cite{ref_Baca} and then also attracted some attention (see e.g. \cite{ref_AnhPal,ref_NurBasSalGao}). The best general result $tvs(G)\leq 3n/\delta$ is due to Anholcer et al. \cite{ref_AnhKalPrz}. Majerski and Przyby{\l}o \cite{ref_Prz4} showed that for dense graphs of sufficiently large order the bound is $tvs(G)\leq (2+o(1))n/\delta+4$.

Let us focus again on the product version of the problem. Anholcer \cite{ref_Anh1} proved several lower bounds for general and regular graphs. In particular, it was shown that for every graph $G$ with $n_d$ vertices of degree $d$,
$$
ps(G)\geq \max_{\delta(G)\leq d\leq\Delta(G)}\left\{\left\lceil \frac{d}{e}n_d^{1/d}-d+1\right\rceil\right\},
$$
while for $r$-regular graphs on $n$ vertices, the inequality takes the form 
$$
ps(G)\geq \left\lceil \frac{r}{e}n^{1/r}-r+1\right\rceil.
$$
In the same paper it was also shown that $ps(G)\leq p(|E(G)|)$, where $p(k)$ denotes the $k^{th}$ prime number. A better bound for sufficiently large graphs follows from the results of Pikhurko \cite{ref_Pikh1}: $ps(G)\leq |E(G)|$. A substantial improvement to this bounds was proved by Darda and Hujdurovi\'c \cite{ref_DarHuj}. Namely, they showed that for most graphs
$$
ps(G) \leq |V(G)| - 1.
$$

In the case of cycles the above results can be improved. In particular, the exact values for short cycles have been presented in \cite{ref_Anh1}. In the same paper it was shown that the lower bound for general cycles is given by the inequality
$$
ps(C_n)\geq  \left\lceil \sqrt{2n}-\frac{1}{2} \right\rceil.
$$

On the other hand, it was proved in \cite{ref_Anh1} that for every $\varepsilon>0$ there exists $n_0$ such that for every $n\geq n_0$
$$
ps(C_n)\leq \lceil (1+\varepsilon)\sqrt{2n}\ln n\rceil.
$$

The two last results hold also for paths $P_n$ and all the Hamiltonian graphs of order $n$. In the same paper it was proved that given any $\varepsilon>0$, one can find $n_j^{(0)}$,$j=1,\dots,k$ such that for every $k$-tuple $(n_1,n_2,\dots,n_k)$,$n_j \geq n_j^{(0)}$, $j=1,2,\dots,k$, the upper bounds for toroidal grids and grids are
$$
ps(T_{n_1 \times n_2 \times \dots \times n_k}) \leq  \lceil (1+\varepsilon)\sqrt{2}(\sum_{j=1}^{k}\sqrt{n_j })\ln{(\sum_{j=1}^{k}n_j)}\rceil
$$
and
$$
ps(G_{n_1 \times n_2 \times \dots \times n_k}) \leq  \lceil (1+\varepsilon)\sqrt{2}(\sum_{j=1}^{k}\sqrt{n_j })\ln{(\sum_{j=1}^{k}n_j)}\rceil,
$$
respectively.

Skowronek-Kazi\'ow considered the local versions of both problems (with vertex labels allowed or not). In both cases the goal was to distinguish only the neighboring vertices. In particular she proved that if the vertex labels are allowed, then $3$ colors are enough to properly color any graph \cite{ref_Skowronek1} and if one cannot label the vertices, then it is enough to use at most $4$ colors \cite{ref_Skowronek2}. Note that the respective lower bounds are $2$ and $3$, so the obtained results are almost optimal. In the case of complete graphs Skowronek-Kazi\'ow obtains these lower bounds, so we can see that
$$
ps(K_n)=3
$$
and
$$
tvps(K_n)=2.
$$

In \cite{ref_Anh2} Anholcer presented, among others, some results for bipartite graphs. In particular, if $2\leq m\leq n$, then
$$
ps(K_{m,n})=3
$$
if and only if $n\leq {{m+2}\choose {2}}$, while otherwise $ps(K_{m,n})\geq 4$. In the same paper  it was shown that given arbitrary integer $D\geq 3$, for almost all forests $F$ such that $\Delta(F)=D$, $n_2=0$, $n_0\leq 1$ and if one removes all the pendant edges, then in the resulting forest $F^\prime$, $n_2=0$, the product irregularity strength equals to
$$
ps(F)=n_1.
$$

Darda and Hujdurovi\'c \cite{ref_DarHuj} generalized some of those results. In particular, they proved that for positive integers $m_1 \leq m_2 \leq m_3 \leq \dots \leq m_k$, and $K_{m_1,m_2,\dots,m_k}$ being the complete multipartite graph such that $m_k\leq m_1+m_2+\dots+m_{k-1}$, the equality $ps(K_{m_1,m_2,\dots,m_k}) = 3$ holds. Baldouski \cite{ref_Bal} in turn showed that for a connected graph $G$ having clique covering number $2\leq \theta(G)\leq 3$, the product irregularity strength $ps(G)=3$, provided that $G$ or the covering cliques have sufficiently large order.

In this paper, some bounds on $tvps(G)$ have been presented. In section \ref{sec_general} we show some general results, as well as the results for regular graphs. In section \ref{sec_cycles} we consider cycles and paths. Then we generalize these results for grids and toroidal grids in section \ref{sec_grids}. Finally some bounds for complete multipartite graphs are presented in section \ref{sec_complete}. We conclude the paper with some open problems.

\section{General Bounds}\label{sec_general}

Since the goal is to obtain distinct product degrees, a unique multiset of labels must be used to label the edges incident to every vertex (although, of course, it does not need to be a sufficient condition). Recall that $n_d$ denotes the number of vertices of degree $d$ (note that $d+1$ labels are present in the product degree of such vertex). Let $s$ be the strength of the labeling (i.e., the masimum label used). Then
$$
n_d\leq{s+d \choose s-1}={s+d \choose d+1} <\left(\frac{e(s+d)}{d+1}\right)^{d+1},
$$
so for every $d$ such that $\delta(G)\leq d\leq\Delta(G)$ the following holds:
$$
s\geq \left\lceil \frac{d+1}{e}n_d^{1/(d+1)}-d\right\rceil.
$$
It means that the following two observations are true.
\begin{myproposition}\label{label_general_lower_bound}
For every graph $G$
$$
ps(G)\geq \max_{\delta(G)\leq d\leq\Delta(G)}\left\{\left\lceil \frac{d+1}{e}n_d^{1/(d+1)}-d\right\rceil\right\}.
$$
\end{myproposition}
\begin{myproposition}
For every $r$-regular graph $G$ such that $|V(G)|=n$,
$$
tvps(G)\geq \left\lceil \frac{r+1}{e}n^{1/(r+1)}-r\right\rceil.
$$
\end{myproposition}

Now let us consider the upper bounds. Obviously, if one assigns the label $1$ to every vertex, then the obtained product degrees are equal to the ones in the case when vertex labels are not allowed. It follows that given any $G$ for which $ps(G)$ is defined, $tvps(G)\leq ps(G)$ holds. In particular, it holds that $ps(G)\leq |V(G)|-1$. However, if one considers a general graph, the best we can obtain is the order of $G$. In fact, after labeling all the edges in an arbitrary way, one can compute the temporary product degrees assuming that all the vertex labels equal to $1$. Now the final labels of vertices, being distinct numbers from the set $\{1,2,\dots,|V(G)|\}$, are assigned in increasing order consistent with non-decreasing order of the temporary product degrees. Obviously the resulting sequence of product degrees is strictly increasing. So we obtain the following.

\begin{myproposition}
For every graph $G$ of order $n$, $tvps(G)\leq n$.
\end{myproposition}

This bound cannot be improved in the general case, as the example of the empty graph shows.

In the next section we provide some results for cycles and paths.

\section{Cycles and Paths}\label{sec_cycles}

Denote by $C_k$ a cycle of length $k$. Then:
\begin{myfact}\label{fact_exact_values}
\begin{eqnarray*}
\begin{array}{l}
tvps(C_{3})=tvps(C_{4})=2,\\
tvps(C_{5})=tvps(C_{6})=tvps(C_{7})=tvps(C_{8})=tvps(C_{9})=tvps(C_{10})=3,\\
tvps(C_{11})=tvps(C_{12})=tvps(C_{13})=tvps(C_{14})=tvps(C_{15})=tvps(C_{16})=4,\\
tvps(C_{n})\geq 5 $ if $n\geq 17.
\end{array}
\end{eqnarray*}
\end{myfact}

\textbf{Proof.} Let us start with proving the fact that the given values of $tvps(C_n)$ are the minimum possible.

By using label $1$, only one product can be obtained: $1$, so it is impossible to label $C_3$ with just this label. Thus, one has to apply at least two labels to label any cycle $C_n$, where $n\geq 3$.

Similarly we can observe that when using labels $1$ and $2$, one can obtain four distinct product degrees: $1$, $2$, $4$ and $8$, so at least three labels must be applied for $C_n$, where $n\geq 5$.

Using labels $1$, $2$ and $3$ it is possible to obtain ten products: $1$, $2$, $3$, $4$, $6$, $8$, $9$, $12$, $18$ and $27$, so we have to use at least four labels to obtain a product-irregular labeling of $C_{11}$ and any longer cycle.

The labels $1$, $2$, $3$ and $4$ can result in at most sixteen product degrees:
$1$, $2$, $3$, $4$, $6$, $8$, $9$, $12$, $16$, $18$, $24$, $27$, $32$, $36$, $48$ and $64$,
so at least $5$ colors are necessary for $C_{n}$, where $n\geq 17$.

Now we are going to show that the values in the theorem are enough to obtain product-irregular labelings. For every cycle $C_n$, $3\leq n \leq 16$, let the sequence of labels be
$$
S_n=([w(v_1)], w(v_1v_2), [w(v2)], \dots, w(v_{n-1}v_n), [w(v_n)], w(v_nv_1))
$$
(the labels of vertices are in the brackets). Sample sequences of labels minimizing the value of $tvps(C_n)$ are listed below:
\begin{eqnarray*}
\begin{array}{l}
S_3=([1],1,[2],2,[2],2),\\
S_4=([1],1,[1],1,[2],2,[2],2),\\
S_5=([1],1,[2],2,[3],3,[3],3,[2],2),\\
S_6=([1],1,[1],1,[2],2,[3],3,[3],3,[2],2),\\
S_7=([1],1,[1],2,[2],2,[2],3,[3],3,[3],2,[2],1),\\
S_8=([1],1,[3],1,[1],2,[2],2,[2],3,[3],3,[3],2,[2],1),\\
S_9=([1],1,[1],3,[2],1,[1],2,[2],2,[2],3,[3],3,[3],2,[2],1),\\
S_{10}=([1],1,[1],3,[1],3,[2],1,[1],2,[2],2,[2],3,[3],3,[3],2,[2],1),\\
S_{11}=([1],3,[1],3,[1],1,[1],1,[1],2,[2],3,[3],3,[3],4,[4],4,[4],3,[3],2),\\
S_{12}=([1],3,[1],3,[1],1,[4],1,[1],1,[1],2,[2],3,[3],3,[3],4,[4],4,[4],3,[3],2),\\
S_{13}=([4],1,[2],1,[2],3,[4],2,[4],2,[3],3,[1],3,[4],4,[2],4,[4],4,[3],3,[3],\\3,[2],2),\\
S_{14}=([4],1,[2],1,[4],1,[2],3,[4],2,[4],2,[3],3,[1],3,[4],4,[2],4,[4],4,[3],\\3,[3],3,[2],2),\\
S_{15}=([4],1,[2],1,[3],1,[4],1,[2],3,[4],2,[4],2,[3],3,[1],3,[4],4,[2],4,[4],\\4,[3],3,[3],3,[2],2),\\
S_{16}=([4],1,[2],1,[1],1,[3],1,[4],1,[2],3,[4],2,[4],2,[3],3,[1],3,[4],4,[2],\\4,[4],4,[3],3,[3],3,[2],2).\\
\end{array}
\end{eqnarray*}
This completes the proof.\qed

Now we are going to show the bounds for arbitrary $n$.

\begin{myproposition}\label{proposition_square}
For every $n>2$,
$$
tvps(C_n)\geq  \left\lceil\sqrt[3]{6n}-1\right\rceil.
$$
\end{myproposition}

\textbf{Proof.}
In the case of $C_n$, the inequality
$$
n_d\leq{s+d \choose d+1}
$$
takes the form
$$
s(s+1)(s+2)\geq 6n,
$$
which implies
$$
(s+1)^3\geq 6n,
$$
and finally
$$
s\geq \sqrt[3]{6n}-1.
$$
\qed

Next, let us prove some upper bounds on $tvps(C_n)$. Recall that in the label sequences, the vertex labels are surrounded with square brackets.

\begin{mylemma}\label{lemma_3s}

For every $s\geq 3$ there exist product-irregular labelings of $C_{3s-2}$, $C_{3s-1}$ and $C_{3s}$ with the label sequence containing the subsequence
$$s-1,[s-1],s,[s],s,[s],s-1$$.
\end{mylemma}
\textbf{Proof.} For $C_7$, $C_8$ and $C_9$, see the sequences given in the proof of Fact \ref{fact_exact_values}. Now, assume that the statement is true for some $s$. Take any $n$, $3s-2\leq n\leq 3s$. Without loss of generality we can assume that $w(v_{n-2}v_{n-1})=s-1$, $w(v_{n-1})=s-1$, $w(v_{n-1}v_{n})=s$, $w(v_{n})=s$, $w(v_{n}v_{1})=s$, $w(v_{1})=s$,$w(v_{1}v_{2})=s-1$. As one can see, the three largest product degrees are $w(v_n-1)=s(s-1)^2$, $w(v_{n})=s^3$ and $w(v_1)=s^2(s-1)$. We now extend the labeling to $C_{n+3}$ adding only one label $s+1$ in the following way (of course the edge $v_nv_1$ is removed together with its label): $w(v_{n}v_{n+1})=s$, $w(v_{n+1})=s$, $w(v_{n+1}v_{n+2})=s+1$, $w(v_{n+2})=s+1$, $w(v_{n+2}v_{n+3})=s+1$, $w(v_{n+3})=s+1$,$w(v_{n+3}v_{1})=s$. As one can easily check, the new labeling contains the sequence
$$s,[s],s+1,[s+1],s+1,[s+1],s.$$
Moreover it preserves the product degrees for all $v_i$, $1\leq i\leq n$, while for the three new vertices we have $pd(v_{n+1})=s^2(s+1)$, $pd(v_{n+2})=(s+1)^3$, $pd(v_{n+3})=s(s+1)^2$. Clearly, the new labeling satisfies the subsequence condition and is product-irregular, so the proof follows by induction.
\qed

\begin{mylemma}\label{lemma_4s}
For every $s\geq 4$ there exist product-irregular labelings of $C_{4s-3}$, $C_{4s-2}$,$C_{4s-1}$ and $C_{4s}$ with the label sequence containing the subsequence
$$s-1,[s],s,[s-2],s,[s],s,[s-1],s-1.$$
\end{mylemma}
\textbf{Proof.} For $C_{13}$, $C_{14}$, $C_{15}$ and $C_{16}$, see the sequences given in the proof of Fact \ref{fact_exact_values}. Now, assume that the statement is true for some $s$. Take any $n$, $4s-3\leq n\leq 4s$. Without loss of generality we can assume that $w(v_{n-2}v_{n-1})=s-1$, $w(v_{n-1})=s$, $w(v_{n-1}v_{n})=s$, $w(v_{n})=s-2$, $w(v_{n}v_{1})=s$, $w(v_{1})=s$,$w(v_{1}v_{2})=s$, $w(v_{2})=s-1$,$w(v_{2}v_{3})=s-1$. As one can see, the four largest product degrees are $w(v_n-1)=s^2(s-1)$, $w(v_{n})=s^2(s-2)$, $w(v_1)=s^3$ and $w(v_2)=s(s-1)^2$. Similarly as in the proof of Lemma \ref{lemma_3s}, we extend the labeling to $C_{n+4}$ adding only one label $s+1$ by assigning $w(v_{n}v_{n+1})=s$, $w(v_{n+1})=s+1$, $w(v_{n+1}v_{n+2})=s+1$, $w(v_{n+2})=s-1$, $w(v_{n+2}v_{n+3})=s+1$, $w(v_{n+3})=s+1$,$w(v_{n+3}v_{n+4})=s+1$, $w(v_{n+4})=s$,$w(v_{n+4}v_{1})=s$. As one can easily check, the new labeling contains the sequence
$$s,[s+1],s+1,[s-1],s+1,[s+1],s+1,[s],s.$$
Moreover it preserves the product degrees for all $v_i$, $1\leq i\leq n$, while for the four new vertices we have $pd(v_{n+1})=s(s+1)^2$, $pd(v_{n+2})=(s-1)(s+1)^2$, $pd(v_{n+3})=(s+1)^3$, $pd(v_{n+4})=s^2(s+1)$. Clearly, the new labeling satisfies the subsequence condition and is product-irregular, so the proof follows by induction.
\qed

The lemmas \ref{lemma_3s} and \ref{lemma_4s} imply the following.
\begin{myproposition}
$$
ps(C_n)\leq
\begin{cases}
\lceil \frac{n}{3}\rceil \text{     if     } n\geq 7\\
\lceil \frac{n}{4}\rceil \text{     if     } n\geq 13
\end{cases}.
$$
\end{myproposition}

The following theorem gives the upper bound for the cycles of sufficiently large order.

\begin{mytheorem}\label{lab_main_theorem}
For every $\varepsilon>0$ there exists $n_0$ such that for every $n\geq n_0$
$$
ps(C_n)\leq \lceil (1+\varepsilon)3\sqrt[3]{2}(1+\varepsilon)n^{1/3}\ln n\rceil.
$$
\end{mytheorem}
\textbf{Proof.} Let $s=\lceil (1+\varepsilon)3\sqrt[3]{2}(1+\varepsilon)n^{1/3}\ln n\rceil$ and let $p$ be the greatest prime that satisfies the condition $p \leq \lfloor 2/3\pi(s)\rfloor$. From Bertrand-Chebyshev theorem (\cite{ref_Tch}), $p \geq \lfloor 2/3\pi(s)\rfloor/2+1$.

We start with finding the edge labels.

For every $q$, $1\leq q<p/2$, we define the sequence: $0$, ${q\mod p}$, ${2q\mod p}$, $\dots$, ${(p-1)q\mod p}$, ${pq\mod p=0}$. We will refer to such a sequence as "the chain".

In any fixed chain, each number $a$, $0\leq a\leq p-1$, appears exactly once, as $q$ and $p$ are relatively prime, and so the order of $q$ in the additive group $Z_p$ is equal to $p-1$. Thus, if we join all those chains together and form a multichain ($0$'s being common members of every two neighboring chains), every pair of numbers from the considered set will appear at most once as a pair of consecutive numbers (here we use the assumption that $q<p/2$). Next, let $p_0=1$, $p_1=2$, $p_2=3$, $p_3=5$ and so on, where the following terms of the sequence are consecutive primes. Replace every number $a$ in the multichain constructed above by $p_a$. None of the pairs repeats, so the same holds for the products. Moreover, since any two consecutive numbers are always distinct, there is no product equal to the square of a natural number. Now we join $r$ identical multichains together in the same way as we joined chains to obtain the multichain and close the cycle, where $1\leq r\leq \lceil\pi(s)\rceil/3$ will be chosen later (in particular it may happen that $r= \lceil\pi(s)\rceil/3$. We are going to show that it is possible to label at least $n$ edges in this way.

In every chain $p$ labels have been used and each chain is based on distinct natural number smaller than $p/2$, so one can label
$$
m\geq r\frac{p(p-1)}{2}\geq \lceil \pi(s)/3\rceil\frac{(2/3\lfloor\pi(s)\rfloor)^2}{8}\geq \frac{1}{54}\left\lfloor\frac{s}{\ln s}\right\rfloor^3
$$
\noindent{}edges. Because of the choice of $s$, for sufficiently large $n$, the inequality $m>n$ holds.

We use the maximum number of multichains and then the maximum number of chains that allows to label at most $n$ edges.  It means that the number of labeled edges is at least $n-(p-1)r/2+1$ but not greater than $n$. Let $t$ denote the number of unlabeled edges. Of course, $0\leq t<(p-1)r/2-1$. Moreover, exactly $p$ labels have been applied and no squares of natural numbers have occured so far. Now let us choose $t$ edges from the last chains of as many multichains as necessary (i.e., choose all the edges in the last chains of $\lfloor t/(p-1)\rfloor$ multichains and $t-(p-1)\lfloor t/(p-1)\rfloor$ in any other mutichain. Then, replace each of those edges with two new edges, assigned the same label. No product degree changes, as the only change is that some sequences of the form $p_i,p_j,p_k$ become $p_i,p_j,p_j,p_k$. Eventually only some new products of the form $p_j^2$ can appear. So we obtain edge labeling of the cycle of length $n$. Note that every product degree appears at most $r$ times, as in every multichain all the degrees are distinct. Now for every set of vertices with equal product degree we use the primes greater than any label used so far (there are exactly $r$ of them) to distinguish the product degrees of those vertices. This results with a product-irregular labeling of $C_n$.
\qed

It is straightforward to see that the above results for cycles $C_n$ can be used also to prove upper bounds on the $tvps(G)$, where $G$ is a path (it is enough to remove any edge with label $1$, Hamiltonian graph (one labels with $1$ all the edges not belonging to the Hamiltonian cycle) or semi-Hamiltonian graph (one labels with $1$ all the edges not belonging to the Hamiltonian path).

\section{Grids and Toroidal Grids}\label{sec_grids}

Given $k$ paths $P_j$ ($j=1,2,\dots k$) with vertex sets $V_j$ ($j=1,2,\dots k$), where $|V_j|=n_j$,the grid $G_{n_1 \times n_2\times \dots \times n_k}$ is the Cartesian product of those paths: $V=V(G_{n_1 \times n_2\times \dots \times n_k})=V_1 \times V_2\times \dots \times V_k$ and two vertices $(u_1, u_2, \dots, u_k)$, $(v_1, v_2, \dots, v_k)\in V$ are adjacent if and only if all their coordinates but one, say $j^{th}$, are the same, and $(u_j,v_j)\in E(P_j)$. By taking cycles instead of paths, one obtains the toroidal grid $T_{n_1 \times n_2\times \dots \times n_k}$.

\begin{mylemma}\label{lab_lemma_sum}
Let $n_1, n_2, \dots, n_k$ be natural numbers, $n_j\geq 3$ for all $j=1, 2, \dots, k$.
\begin{enumerate}

\item
Let $p_1, p_2, \dots, p_k$ and $r_1, r_2, \dots, r_k$ be $2k$ natural numbers (not necessarily distinct) such that $p_j$ primes as the edge labels and $r_j$ primes as vertex labels are enough to obtain a product-irregular labeling of $C_{n_j}$ ($j=1,2,\dots k$). Then $\sum_{j=1}^{k}{(p_j+r_j)}$ primes are enough to obtain a product-irregular labeling of the toroidal grid $T_{n_1 \times n_2\times \dots \times n_k}$.

\item
Let $p_1, p_2, \dots, p_k$ and $r_1, r_2, \dots, r_k$ be $2k$ natural numbers (not necessarily distinct) such that $p_j$ primes as the edge labels and $r_j$ primes as vertex labels are enough to obtain a product-irregular labeling of $P_{n_j}$ ($j=1,2,\dots k$). Then $\sum_{j=1}^{k}{(p_j+r_j)}$ primes are enough to obtain a product-irregular labeling of the grid $G_{n_1 \times n_2\times \dots \times n_k}$.

\end{enumerate}

\end{mylemma}

\textbf{Proof.} First consider the toroidal grids. We will use the induction on $k$. For $k=1$ it is trivially true. Assume that we have managed to label $T_{n_1 \times n_2\times \dots \times n_{k-1}}$ using $\sum_{j=1}^{k-1}{p_j}$ primes. To construct $T_{n_1 \times n_2\times \dots \times n_k}$, join $n_k$ copies of $T_{n_1 \times n_2\times \dots \times n_{k-1}}$ in such a way that for every vertex $v$ of $T_{n_1 \times n_2\times \dots \times n_{k-1}}$, all copies of $v$ are joined with the same copy of $C_{n_k}$. The sets of labels applied to the edges incident with each copy of $v$ are distinct, as they consist of constant set of labels used to label $T_{n_1 \times n_2\times \dots \times n_{k-1}}$ and distinct pairs of labels used on $C_{n_k}$ (taken from a disjoint set of numbers). We define the vertex label as the product of the labels of the corresponding vertices in the factors (i.e., in $T_{n_1 \times n_2\times \dots \times n_{k-1}}$ and $C_{n_k}$). It results in a product-irregular labeling of $T_{n_1 \times n_2\times \dots \times n_k}$ using $p_k+r_k$ new primes, which is the desired number of labels.

In the case of $G_{n_1 \times n_2\times \dots \times n_k}$ we use $P_{n_j}$ instead of $C_{n_j}$.\qed

Note that the maximum label used in the above labeling is not greater than $\max\{P,R\}$, where $P$ is the maximum edge label among all cycles (paths) and $R$ the product of the maximum vertex labels of all cycles (paths). Now we can formulate the main result of this section.

\begin{mytheorem}
For every $k\geq 2$ and every $\varepsilon>0$ there exist $n_j^{(0)}$,$j=1,\dots,k$ such that for every $k$-tuple $(n_1,n_2,\dots,n_k)$,$n_j \geq n_j^{(0)}$, $j=1,2,\dots,k$,
\begin{enumerate}

\item
$ps(T_{n_1 \times n_2 \times \dots \times n_k}) \leq  \lceil (1+\varepsilon)k^{k-2}(\sum_{j=1}^{k}n_j^{k/(2k+1)})\sum_{j=1}^{k}\ln{n_j}\rceil$;

\item
$ps(G_{n_1 \times n_2 \times \dots \times n_k}) \leq  \lceil (1+\varepsilon)k^{k-2}(\sum_{j=1}^{k}n_j^{k/(2k+1)})\sum_{j=1}^{k}\ln{n_j}\rceil$.

\end{enumerate}
\end{mytheorem}

\textbf{Proof.} In order to find an irregular labeling of a cycle $C_n^j$ it is enough to use $p_j$ labels (primes or $1$) to label the edges and
$$
r_j=\left\lceil\frac{2n}{p_j(p_j-1)}\right\rceil
$$
labels for vertices (see the proof of Theorem \ref{lab_main_theorem}). In particular, if $p_j= \lfloor n_j^{k/(2k+1)}\rfloor$, then $r_j\leq \lceil n_j^{1/(2k+1)}\rceil$. This means that we need at most
$$
\sum_{j=1}^{k}{(\lfloor n_j^{k/(2k+1)}\rfloor+\lceil n_j^{1/(2k+1)}\rceil)}
$$
distinct labels (all of them being primes or $1$). Let
$$
s=\lceil (1+\varepsilon)k^{k-2}(\sum_{j=1}^{k}n_j^{k/(2k+1)})\sum_{j=1}^{k}\ln{n_j}\rceil
$$
and
$$
\varepsilon^\prime=(1+\varepsilon)^{1/k}-1.
$$
We assign
$$
\sum_{j=1}^{k}\left\lceil n_j^{1/(2k+1)}\right\rceil
$$
smallest primes to the vertices ($1$ will be used as edge label). This means that the maximum vertex label $R_j$ in every cycle satisfies the condition
$$
R_j<(1+\varepsilon^\prime)\left(\sum_{j=1}^{k}n_j^{1/(2k+1)}\right)\ln{\left(\sum_{j=1}^{k}n_j^{1/(2k+1)}\right)}
$$
for sufficiently large $s_j$ (i.e. $n_j$), $j=1,\dots,n$. Consequently, the maximum vertex label $R$ in $T_{n_1 \times n_2 \times \dots \times n_k}$ is not greater than
\begin{align*}
R=\prod_{j=1}^{k}{R_j}< (1+\varepsilon)\left(\left(\sum_{j=1}^{k}n_j^{1/(2k+1)}\right)\ln{\left(\sum_{j=1}^{k}n_j^{1/(2k+1)}\right)}\right)^k
\end{align*}
for sufficiently large $n_j$, $j=1,\dots, k$. For sufficiently large $n_j$ we have
$$
\ln{\left(\sum_{j=1}^{k}n_j^{1/(2k+1)}\right)}\leq \sum_{j=1}^{k}\ln{\left(n_j^{1/(2k+1)}\right)}.
$$
From H\"older's Inequality it follows that
$$
\sum_{j=1}^{k}{n_j^{1/{(2k+1)}}}\leq \left(\sum_{j=1}^{k}{n_j^{k/{(2k+1)}}} \right)^{1/k}k^{1-1/k}
$$
and
\begin{align*}
\sum_{j=1}^{k}\ln{\left(n_j^{1/(2k+1)}\right)}\leq \left(\sum_{j=1}^{k}\left(\ln{\left(n_j^{1/(2k+1)}\right)}\right)^k\right)^{1/k}k^{1-1/k}\\
<\left(\sum_{j=1}^{k}\ln{\left(kn_j^{1/(2k+1)}\right)}\right)^{1/k}k^{1-1/k}\\
<\left(\sum_{j=1}^{k}\ln{\left(n_j^{k/(2k+1)}\right)}\right)^{1/k}k^{1-1/k}
\end{align*}
This implies
$$
R<(1+\varepsilon)k^{2k-2}\left(\sum_{j=1}^{k}{n_j^{k/{(2k+1)}}} \right) \sum_{j=1}^{k}\ln{n_j}\leq s.
$$
On the other hand, for sufficiently large $n_j$, $j=1,\dots, n$, the highest edge label is not greater than
\begin{align*}
(1+\varepsilon)\left(\sum_{j=1}^{k}{(\lfloor n_j^{k/(2k+1)}\rfloor+\lceil n_j^{1/(2k+1)}\rceil)}\right)\ln{\left(\sum_{j=1}^{k}{(\lfloor n_j^{k/(2k+1)}\rfloor+\lceil n_j^{1/(2k+1)}\rceil)}\right)}\\
<(1+\varepsilon)\left(2\sum_{j=1}^{k}{ n_j^{k/(2k+1)}}\right)2\sum_{j=1}^{k}{\ln{n_j^{k/(2k+1)}}}<s.
\end{align*}
Since in every cycle $1$ appears as an edge label, the same reasoning proves the theorem for the grid $G_{n_1 \times n_2\times \dots \times n_k}$.
\qed

\section{Complete multipartite graphs}\label{sec_complete}

In this section we prove two simple results for complete multipartite graphs. Let us start with the bipartite case.
\begin{myproposition}
Let m and n be two integers such that $3\leq m\leq n\leq \binom{m+2}{2}$. Then\\
$$tvps(K_{m,n})=3.$$
\end{myproposition}
\textbf{Proof.} The fact that the labels $1$, $2$ and $3$ are enough follows from the equality $ps(K_{m,n})=3$ (it is enough to label all vertices with $1$ and we are done). On the other hand, if we use only labels $1$ and $2$, then the product degrees could take at most $n+2$ values, namely $1, 2, 4, \dots, 2^{n+1}$, while there are at least $n+3$ vertices, so at least third label is necessary.
\qed.

By using a similar argument, we can prove the exact value for complete multipartite graphs. 

\begin{myproposition}
Let $m_1,m_2,\dots, m_k$ be integers such that $3\leq m_1 \leq m_2 \leq m_3 \leq \dots \leq m_k \leq m_1+m_2+\dots+m_{k-1}$. Then\\
$$tvps(K_{m_1,m_2,\dots,m_k}) = 3.$$
\end{myproposition}
\textbf{Proof.} Let us recall that Darda and Hujdurovi\'c \cite{ref_DarHuj} proved that if $K_{m_1,m_2,\dots,m_k}$ is the complete multipartite graph such that $0<m_1 \leq m_2 \leq m_3 \leq \dots \leq m_k \leq m_1+m_2+\dots+m_{k-1}$, then the equality $ps(K_{m_1,m_2,\dots,m_k}) = 3$ holds. By putting label $1$ on every vertex we get $tvps(K_{m_1,m_2,\dots,m_k}) \leq 3$. In order to show that $2$ labels are not enough, observe that $\Delta(K_{m_1,m_2,\dots,m_k}) = m_2+\dots+m_k$, so using labels $1$ and $2$ one can obtain at most $2+m_2+\dots+m_k$ different weights, while the number of vertices is $m_1+m_2+\dots+m_k\geq 3+m_2+\dots+m_k$.
\qed.

\section{Conclusion}

In this paper, a new graph invariant, total vertex product irregularity strength, has been introduced. It is a parameter similar to the total vertex irregularity strength, where where the weighted degrees are computed as products. We proved several results for general and regular graphs, as well as for some specific families of graphs, like cycles, paths and grids.

Most of our results show that for any $d$-regular graph $G$ of order $n$ it holds that $c_1n^{1/(d+1)}\leq g\leq c_2n^{1/(d+1)}\ln{n}$ for some constants $c_1$ and $c_2$. Erd\"os \cite{ref_Erd} proved that in the case of integers not greater than $n$, the cardinality $s(n)$ of a subset resulting in distinct pairwise products (multiplicative Sidon set) cannot be much greater than the number of primes not greater than $n$:
$$
s(n)=\pi(n)+\Theta\left(\frac{n^{3/4}}{\ln^{3/2}n}\right).
$$
This makes us believe that the upper bounds presented in this paper are closer to the exact value of $tvps(G)$ than the lower and the effort of the further research should be focused on the improvement of the latter ones. For that reason we formulate the two following open problems.
\begin{myproblem}
Is there a constant $c$ such that for every $d$-regular graph $G$ of order $n$
$$
tvps(G)\geq cn\ln{\ln{n}}?
$$
\end{myproblem}
\begin{myproblem}
Is there a constant $c$ such that for every $d$-regular graph $G$ of order $n$
$$
tvps(G)\geq \frac{cn\ln{n}}{\ln{\ln{n}}}?
$$
\end{myproblem}
We also believe that the logarithmic factor is enough to guarantee the existence of a product-irregular labeling of any graph. This encourages us to pose the following conjecture.
\begin{myconjecture}
There exists a constant $c$ such that for every $d$-regular graph $G$ of order $n$
$$
tvps(G)\leq cn^{1/(d+1)}\ln{n}.
$$
\end{myconjecture}

\section*{Acknowledgment}

This work was partially done while the second author (A. S. Emadi) was visiting the Pozna\'n University of Economics and Business, supported by the University of Mazandaran and the Ministry of Science, Research and Technology of Iran, under the program Ph.D. Research Opportunity.

\end {document}